\documentstyle[11pt]{article}
\newlength{\standardunitlength}
 \newtheorem{lemma}{Lemma}
\newtheorem{theorem}{Theorem} 
\newenvironment{proof}{\noindent {\sc Proof:}}{$\Box$ \vspace{2 ex}}

\begin{document}

\begin{center}
A Generalization of an Alternating Sum Formula for Finite Coxeter Groups
\end{center}

\begin{center}
By Jason Fulman
\end{center}

\begin{center}
Dartmouth College
\end{center}

\begin{center}
Department of Mathematics
\end{center}

\begin{center}
6188 Bradley Hall
\end{center}

\begin{center}
Hanover, NH 03755, USA
\end{center}

\begin{center}
email:jason.e.fulman@dartmouth.edu
\end{center}

\begin{center}
\end{center}

\begin{abstract}
	For $W$ a finite Coxeter group, a formula is found for the size
of $W$ equivalence classes of subsets of a base. The proof is a case-by-case
analysis using results and tables of Orlik and Solomon. As a corollary we obtain
an alternating sum identity which generalizes a well-known identity from the
theory of Coxeter groups.
\end{abstract}

\section{Introduction}

	Let $W$ be a finite Coxeter group and $\Pi$ a set of fundamental
roots or base of $W$. Note that $W$ acts on subsets of $\Pi$. For $J,K \subseteq
\Pi$, one writes $J \sim K$ if $J=w(K)$ for some $w \in W$. Letting $\lambda(K)$
denote the equivalence class of $K$ under this action, a natural problem is to
compute the size $|\lambda(K)|$ of this equivalence class. It suffices to work
with irreducible $W$. For the classical types $|\lambda(K)|$ can be computed
without much difficultly; for the exceptional types $|\lambda(K)|$ has been
essentially computed by Carter \cite{C}.

	Theorem \ref{equality} of this note will give a formula for $|\lambda(K)|$ in
terms of the index of the parabolic subgroup $W_K$ in its normalizer $N_W(W_K)$
and the characteristic polynomial $\chi({\cal L}^K,t)$ of an upper interval in the
intersection lattice $L$ of root hyperplanes corresponding to $W$. This
result is of interest in that it gives an expression for $|\lambda(K)|$ in terms
of quantities naturally associated to $W$ and $K$. Our proof is a case-by-case
analysis; a unified proof for all finite Coxeter groups would be desirable.

	Theorem \ref{generalize} applies Theorem \ref{equality} to obtain a
generalization of the alternating sum formula

\[ \sum_{K \subseteq \Pi} (-1)^{|K|} \frac{|W|}{|W_K|} = 1. \]

	The above alternating sum formula is discussed in Sections 1.15-1.16 of
Humphreys \cite{H}. Topological proofs using the Hopf trace formula are
given by Solomon \cite{So} and Steinberg \cite{St}. Chapter 3 of Humphreys
gives applications to the invariant theory of Coxeter groups, and Section 6.2 of
Carter \cite{C2} gives an application to proving the irreducibility of the
Steinberg character. We hope that our generalization will have similar
interpretations and applications.

\section{Notation} 

 We collect some notation which will be indispensable for what follows. Denote
the set of root hyperplanes of $W$ by $\cal{A}$. Let $L$ be the set of
intersections of the hyperplanes in $\cal{A}$, taking $V
\in  L$. Partially order $L$ by reverse inclusion. Recall that the Moebius
function $\mu$ is defined by $\mu(X,X)=1$ and $\sum_{X \leq Z
\leq Y} \mu(Z,Y)=0$ if $X<Y$ and $\mu(X,Y)=0$ otherwise. The characteristic
polynomial of $L$ is defined as

	\[ \chi(L,x) = \sum_{X \in L} \mu(V,X) x^{dim(X)}. \]

	For $K \subseteq \Pi$, let $Fix(W_K)$ be the fixed space of the parabolic
subgroup $W_K$ in its action on $V$. Let $L^{Fix(W_K)}$ denote the poset
isomorphic to the segment $\{Y \in L(\cal{A})$$ |Y \geq Fix(W_K)\}$.

	Let $n$ be the rank of $W$, and let $N_W(W_K)$ be the normalizer of $W_K$ in
$W$. As in the introduction, for $J,K \subseteq \Pi$ write $J \sim K$ if
$J=w(K)$ for some $w \in W$ and let $\lambda(K)$ denote the equivalence class
of $K$ under this action.

\section{Main Results}

	The following lemma of Bergeron, Bergeron, Howlett, and Taylor translates the
equivalence relation $\sim$ into a condition about conjugacy of parabolic
subgroups.

\begin{lemma} (Bergeron, Bergeron, Howlett, Taylor \cite{BBHT}) If $J,K \subseteq
\Pi$, then $J \sim K$ if and only if $W_J$ and $W_K$ are conjugate.
\end{lemma}

	Now we derive an expression for $|\lambda(K)|$ in terms of quantities naturally
associated to $W$.	

\begin{theorem} \label{equality} Let $W$ be a finite Coxeter group of rank $n$.
Then for all $K \subseteq \Pi$,

\[ |\lambda(K)| = (-1)^{n-|K|} \frac{|W_K|}{|N_W(W_K)|} \chi(L^{Fix(W_K)},-1). \]

\end{theorem}

\begin{proof}
	As both sides of the conjectured equation are multiplicative with respect
to direct product of groups, it suffices to prove the theorem for irreducible
$W$. This will be done case by case.

	For $W$ of type $A_{n-1}$, it is
easy to see that $J \sim K$ exactly when $W_J$ and $W_K$ are isomorphic.
Given $K$, define equivalence classes called blocks on the set $\{1,\cdots,n\}$ by
letting $p \sim q$ if it is possible to transpose $p$ and $q$ using only elements
in $W_K$. Let $n_i$ be the number of blocks of size $i$. A subset
$J \subseteq \Pi$ such that $W_J$ is isomorphic to $W_K$ arises from any of the
$(n- \sum in_i)!$ permutations of the blocks of $K$. One must then divide by
$\prod_i n_i!$ since permuting the blocks of size
$i$ amongst themselves leads to the same $J$. Thus,

\begin{eqnarray*}
|\lambda(K)| & = & \frac{(n-\sum_i n_i)!}{\prod_i n_i!}\\
& = & \frac{(n-|K|)!}{\prod_i n_i!}.
\end{eqnarray*}

	Since $N_W(W_K)$ is isomorphic to $\prod_i [(A_i)^{n_i} \times A_{n_i}]$, one
has that $\frac{|W_K|}{|N_W(W_K)|} = \frac{1}{\prod_i n_i!}$. Proposition 2.1 of
Orlik and Solomon \cite{OS} states that $\chi(L^{Fix(W_K)},t) = (t-1)
\cdots (t-(n-|K|-1))$. Thus,

\[ (-1)^{n-|K|} \frac{|W_K|}{|N_W(W_K)|} \chi(L^{Fix(W_K)},-1) = \frac{(n-|K|)!}{\prod_i
n_i!}. \]

	Next we prove the result for type $B_n$ (the proof for type $C_n$ is identical).
It is readily checked that $K \sim J$ exactly when $W_K$ and $W_J$ are isomorphic.
Let $W_K$ be isomorphic to $B_j \times A_{\lambda_1-1} \times A_{\lambda_2-1}
\times \cdots$ and let $n_i$ be the number of parts of $\lambda$ of size $i$.
Arguing as for type $A$,

\begin{eqnarray*}
|\lambda(K)| & = & \frac{(n-j-\sum in_i)!}{\prod_i n_i!}\\
& = & \frac{(n-|K|)!}{\prod_i n_i!}.
\end{eqnarray*}

 On page 12 of Carter \cite{CProc}, it is proved that $\frac{|W_K|}{|N_W(W_K)|}
= \frac{1}{\prod_i 2^{n_i} n_i!}$. Proposition 2.2 of Orlik and Solomon \cite{OS}
states that $\chi(L^{Fix(W_K)},t) = (t-1)(t-3) \cdots (t-(2(n-|K|)-1)$. Therefore
we conclude as desired that

\[ (-1)^{n-|K|} \frac{|W_K|}{|N_W(W_K)|} \chi(L^{Fix(W_K)},-1) = \frac{(n-|K|)!}{\prod_i
n_i!}. \]

	Next we proceed to type $D_n$. Take $\Pi = \{\alpha_1,\cdots,\alpha_n\}$ a
fundamental set of simple roots with $\alpha_1=\varepsilon_1-\varepsilon_2,
\alpha_2=\varepsilon_2-\varepsilon_3, \cdots , \alpha_{n-1}=\varepsilon_{n-1} -
\varepsilon_n, \alpha_n = \varepsilon_{n-1}+\varepsilon_n$.

	There are two subcases. The first subcase is that $\{ \alpha_{n-1},\alpha_n \}
\subseteq K$. Let $W_K$ be isomorphic to $D_j \times A_{\lambda_1-1} \times
A_{\lambda_2-1} \times \cdots$ (here we insist that $j+\sum_i \lambda_i=n$) and
let $n_i$ be the number of parts of $\lambda$ of size $i$. Arguing as in type $A$
shows that $|\lambda(K)| = \frac{2(n-|K|)!}{\prod_i n_i!}$. The extra factor of
$2$ comes from the symmetry of the Dynkin diagram which switches $\alpha_{n-1}$
and $\alpha_n$ but leaves the other roots fixed.

	On page 11 of Carter \cite{CProc}, it is proved that $\frac{|W_K|}{|N_W(W_K)|}
= \frac{2}{\prod_i 2^{n_i} n_i!}$. Proposition 2.6 of Orlik and Solomon \cite{OS}
states that $\chi(L^{Fix(W_K)},t) = (t-1)(t-3) \cdots (t-(2(n-|K|)-1)$. Thus

\[ (-1)^{n-|K|} \frac{|W_K|}{|N_W(W_K)|} \chi(L^{Fix(W_K)},-1) = \frac{2(n-|K|)!}{\prod_i
n_i!}. \]

	The second subcase is that $\{ \alpha_{n-1},\alpha_n \} \not \subseteq K$. All
$J$ such that $J \sim K$ also satisfy $\{ \alpha_{n-1},\alpha_n \} \not \subseteq
J$. Given $K$, define equivalence classes called blocks on the set
$\{1,\cdots,n\}$ by letting $p \sim q$ if it is possible to transpose $p$ and $q$
up to sign change, using only elements in $W_K$. Let $n_i$ be the number of
blocks of size $i$. Observe that

\begin{eqnarray*}
|\lambda(K)| & = & |\{J: J \sim K, J \subseteq \{\alpha_1,\cdots,\alpha_{n-1}\}| +
|\{J: J \sim K, J \subseteq \{\alpha_1,\cdots,\alpha_{n-2},\alpha_n\}|\\
& & - |\{J: J \sim K, J \subseteq \{\alpha_1,\cdots,\alpha_{n-2}\}|\\
& = & 2 \frac{(n-|K|)!}{\prod_{i \geq 1} n_i!} -
\frac{(n-|K|-1)!} {(n_1-1)! \prod_{i \geq 2} n_i!}\\
& = & \frac{[2(n-|K|)-n_1] (n-|K|-1)!}{\prod_{i \geq 1} n_i!}.
\end{eqnarray*}

	On page 11 of Carter \cite{CProc}, it is proved that $\frac{|W_K|}{|N_W(W_K)|}
= \frac{2}{\prod_i 2^{n_i} n_i!}$. Proposition 2.6 of Orlik and Solomon \cite{OS}
states that $\chi(L^{Fix(W_K)},t) = (t-1)(t-3) \cdots (t-(2(n-|K|)-3)(t-(n-|K|+r-1)$,
where $r$ is the number of blocks of size greater than one. Clearly
$r=n-|K|-n_1$. Therefore,

\begin{eqnarray*}
(-1)^{n-|K|} \frac{|W_K|}{|N_W(W_K)|} \chi(L^{Fix(W_K)},-1) & = & \frac{2}{\prod_i
2^{n_i} n_i!} [(n-|K|)-n_1] (n-|K|-1)!\\ & = & |\lambda(K)|.
\end{eqnarray*}

	Next we proceed to the cases $I_2(m)$ and $G_2$. On page 175 of Orlik and
Solomon \cite{OS}, $\chi(L^{Fix(W_K)},t)$ is computed in the following relevant
cases:

\begin{itemize}

\item If $n-|K|=0$ then $\chi(L^{Fix(W_K)},t)=1$.
\item If $n-|K|=1$ then $\chi(L^{Fix(W_K)},t)=t-1$.
\item If $n-|K|=2$ then $\chi(L^{Fix(W_K)},t)=(t-1)(t-n)$, where $n+1$
is the number of lines of $L$ contained in $Fix(W_K)$.

\end{itemize}

	As an application of this computation, we prove the theorem for $I_2(m)$ with
$m$ even ($I_2(m)$ with $m$ odd and $G_2$ are easy and similar). One checks that
$|\lambda(K)|=1$ for all $K$ contained in the two element set $\Pi=
\{\alpha_1,\alpha_2\}$. For $K = \emptyset$, $\frac{|W_K|}
{|N_W(W_K)|}=\frac{1}{2m}$ and $Fix(W_K)$ contains $m$ lines of $L$, so the
theorem checks. For $K$ such that $|K|=1$, $\frac{|W_K|}{|N_W(W_K)|}=\frac{1}{2}$
and $\chi(L^{Fix(W_K)},t)=t-1$, so the theorem checks. Finally, for $K=\Pi$,
$\frac{|W_K|}{|N_W(W_K)|}$ and $\chi(L^{Fix(W_K)},t)$ are both equal to $1$.

	The proof of the theorem for the exceptional cases $F_4,H_3,H_4,E_6,E_7,E_8$
consists of a finite number of straightforward calculations based on tables 3-8 of
Orlik and Solomon \cite{OS}. We will work out the details for an example to
illustrate what is involved.

	For the example, take $W=E_7$ and $K$ any subset of $\Pi$ whose Dynkin
diagram has type $(A_1)^2$. As the first column of their Table 7 indicates, all
such $K$ are equivalent under $\sim$. A glance at the Dynkin diagram of $E_7$ thus
shows that $|\lambda(K)|=15$. The third entry in the first row of their Table 7
gives that $\frac{|W|}{|N_W(W_K)|}=945$. The end of the third row in their Table 7
gives that $\chi(L^{Fix(W_K)},t)=(t-1)(t-5)(t-7)(t-9)(t-11)$. This data shows
that the theorem checks.

	The only minor complications arise in computing $|\lambda(K)|$ in cases where
there are subsets $J$ with Dynkin diagram isomorphic to $K$, but such that $J
\not \sim K$. As the tables of Orlik and Solomon (loc. cit.) indicate, this
happens only for $W=E_7$, and the values $|\lambda(K)|$ for these cases appear on
the top of page 279 of their article.
\end{proof}
	
	From its statement, it is not evident that Theorem \ref{equality} generalizes any
known identities about Coxeter groups. Theorem \ref{generalize} shows that this is
indeed the case.

\begin{theorem} \label{generalize} Let $W$ be a finite Coxeter group
of rank $n$ with base $\Pi$. Then

\[ \sum_{K \subseteq \Pi} (-1)^{n-|K|} \frac{|W|}{|W_K|}
\frac{\chi(L^{Fix(W_K)},t)}{\chi(L^{Fix(W_K)},-1)} = t^n. \]

\end{theorem}

\begin{proof}
	As is explained on page 274 of Orlik and Solomon \cite{OS}, elementary
properties of Moebius functions yield the identity

\[ \sum_{Y \in L} \chi(L^Y,t) = t^n. \]

	Observe that $W$ acts on the lattice $L$. Let ${\cal O}(Fix(W_K))$ be
the orbit of $Fix(W_K)$ under this action. Lemma 3.4 of Orlik and Solomon (loc.
cit.) states that $|{\cal O}(Fix(W_K))| = \frac{|W|}{|N_W(W_K)|}$. For
$Y \in L(A)$, Lemma 3.4 of Orlik and Solomon (loc. cit.) shows that $Fix(W_Y)=Y$.
Since the parabolic subgroup $W_Y$ is conjugate to a standard parabolic subgroup
$W_K$ for some $K$, it follows that $Y \in {\cal O}(Fix(W_K))$ for some $K$.
Then $Y \in {\cal O}(Fix(W_J))$ for exactly the $|\lambda(K)|$ many $J$'s such
that $J \sim K$. Therefore, the Orlik-Solomon identity becomes

\[ \sum_{K \subseteq \Pi} \frac{|W|}{|N_W(W_K)|} \frac{\chi(L^{Fix(W_K)},t)}
{|\lambda(K)|} = t^n.\]

	The result follows from the formula for $|\lambda(K)|$ in Theorem \ref{equality}.
\end{proof} 

{\bf Remarks:}

\begin{itemize}

\item Observe that setting $t=-1$ in Theorem \ref{generalize} yields the
well-known identity

\[ \sum_{K \subseteq \Pi} (-1)^{|K|} \frac{|W|}{|W_K|} = 1. \]

	As noted in the introduction, this identity has a topological proof and arises in
the invariant theory of $W$. It would be desirable to find
analogous interpretations for Theorem \ref{generalize}.

\item The identity

\[ \sum_{K \subseteq \Pi} (-1)^{|K|} \frac{|W|}{|W_K|} = 1 \]

has two known generalizations. To describe the first, let $d_i(W)$ be the
degrees of a Coxeter group $W$. Then from Sections 1.11 and 3.15 of Humphreys
\cite{H},

\[ \sum_{K \subseteq \Pi} (-1)^{|K|} \prod_{i=1}^n
\frac{t^{d_i(W)}-1}{t-1} / \prod_{i=1}^{|K|} \frac{t^{d_i(W_K)}-1}{t-1} = t^n. \]

	To describe the second, note that $W$ acts on the left cosets $vW_K$ by left
multiplication. Let $f_K(w)$ be the number of left cosets of $W_K$ fixed by $w$.
Let $det(w)$ be the determinant of $w$ in its action on $V$. Proposition 1.16 of
Humphreys \cite{H} states that

\[ \sum_{K \subseteq \Pi} (-1)^{|K|} f_K(w) = det(w). \]

	What is the relation of Theorem \ref{generalize} with these generalizations?
Can Theorem \ref{generalize} be further extended to include them?

\end{itemize}

\end{document}